\documentstyle{amsppt}
\NoBlackBoxes
\ifx\undefined\rom
  \define\rom#1{{\rm #1}}
\fi
\ifx\undefined\curraddr
  \def\curraddr#1\endcurraddr{\address {\it Current address\/}: #1\endaddress}
\fi
 

\topmatter
\title Musings on Magnus
\endtitle
\author Gilbert Baumslag
\endauthor
\address Department of Mathematics, City College of New York,
New York, N.Y. 10031
\endaddress 
\keywords Residually torsion-free nilpotent groups, one-relator groups, $\Cal D$-groups\endkeywords
\subjclass Primary 20F05, 20F14 \endsubjclass
\abstract  The object of this paper is to describe a simple
method for proving that certain groups are residually torsion-free nilpotent,
to describe some new parafree groups and to raise some new problems
in honour of the memory of Wilhelm Magnus.
\endabstract
 \thanks The  author was supported in part by  NSF
Grant \#9103098
\endthanks
\endtopmatter
 
\document
\head 1. Introduction
\endhead
I first heard of Wilhelm Magnus  in 1956, when I was
attending some lectures by B.H. Neumann on amalgamated products. At some
point during the course of these lectures, Neumann
remarked that Magnus was the first mathematician to recognize the
value of amalgamated products and had shown just how
effective   a tool they were, in his work on groups defined by a single relation. 
 I was working, at that time, 
 on an universal algebra variation of  free groups, involving
groups with unique roots, which I called $\Cal D$-groups 
 \cite{1}. 
Consequently,   
in an attempt to find analogues of various theorems about 
free groups for $\Cal D$-groups,
I found myself
reading a  beautiful paper of Magnus, in which he
proved the residual torsion-free nilpotence of free groups.
Much of my talk
today will be concerned with these two topics, one-relator groups and
residual nilpotence.
Let me begin by
reminding  you of some of the definitions involved. 
Let $\Cal P$ be a property of groups. 
\definition{Definition} We say that a group $G$ 
is residually a $\Cal P$-group
if for each $g \in G$, $g \neq 1$, there exists a normal subgroup N of
G, such that $g \notin N$ and $G/N$ has $\Cal P$.
\enddefinition 
 The properties that I will be mainly concerned with here 
are {\it freeness}, {\it nilpotence} and
{\it torsion-free nilpotence}. I will make use of the usual commutator notation.
 Thus if
H and K are subgroups of a group G then
$$[H,K]= gp(h^{-1}k^{-1}hk \mid h \in H, k \in K)$$
is the subgroup generated by all the commutators $h^{-1}k^{-1}hk$.
The {\it lower central series of G} is defined to be the series
$$G=\gamma_1(G)\geq \gamma_2(G) \geq \dots \geq \gamma_n(G) \geq \dots,$$
where $\gamma_{n+1}(G)=[\gamma_n(G),G]$. G is termed {\it nilpotent} if
$\gamma_{c+1}(G)=1$ for some $c$. 
 I am now in a position to formulate the  theorem of Magnus \cite{12} 
that I alluded to before.
\proclaim{Theorem 1}
Free groups are residually torsion-free nilpotent.
\endproclaim
The basic idea involved in the proof of Theorem 1 is beautifully simple.
Magnus concocts a faithful representation of a given
free group $F$ in the group of units of a carefully
chosen ring $R$ with 1.  Each
of the elements  $f \in F$ takes the form $1+\phi$, where
$\phi$ lies in an ideal $R^+$ of $R$. $R^+$ carries 
with it the structure of a metric space designed so as to ensure
that if $f \in \gamma_n(F)$, then $d(\phi,0) \leq 2^{-n}$, where
here $d(\phi,0)$ denotes the distance between $\phi$ and $0$.
This suffices to ensure that
$$\bigcap_{n=1}^{\infty} \gamma_n(F)=1.$$
The sketch of the proof,  below, amplifies these remarks.
\demo{Proof} Let F be free on $x_1,\dots,x_q$. Consider the ring $R$ of
power series in the non-commuting indeterminates $\xi_1,\dots,\xi_q$
with rational coefficients. 
Each element $r \in R$ can be thought of as an infinite sum 
which takes the form
$$r=r_0+r_1+\dots+r_n+\dots ,$$
where
$r_0\in \Bbb Q$  
and the {\it homogeneous component} $r_n$ of $r$ of degree $n>0$
 is a finite sum of rational multiples
$c\xi_{i_1}\dots \xi_{i_n}$ of {\it  monomials} $\xi_{i_1}\dots \xi_{i_n}$
of degree $n$, $c \in \Bbb Q, i_j\in \{1,\dots,q\}$. Here 
$\xi_{i_1}\dots \xi_{i_n}=\xi_{j_1}\dots \xi_{j_m}$
only if $n=m$ and $i_r=j_r$ for each $r=1,\dots,n$.
Each of the elements $a_i=1+\xi_i$ is invertible
in $R$, with inverse $a_i^{-1}=1-\xi_i+\xi_i^2-\xi_i^3 + \dots$.
 It then turns out that the elements $a_1,\dots,a_q$ 
 freely generate
a free subgroup of rank $n$  of the multiplicative group of
units of $R$. We define $R^+$ to be the ideal of $R$ consisting of
those elements $r \in R$ such that $r_0=0$. 
We define a metric  $d$ on $R^+$ by setting $d(r,0)=0$ if $r=0$ and,
if $r \ne 0$,
$$d(r,0)= 2^{-n}$$
where $n$ is the degree of
 the first non-zero homogeneous component of $r$. Magnus  proves 
 that if $f = 1+\phi \in \gamma_n(F)$,
then all of the homogeneous components $f_1=f_2=\dots=f_{n-1}=0$, i.e.,
$d(\phi,0)\leq 2^{-n}$.
It is not hard then to deduce that  $F$ is residually torsion-free nilpotent.
\enddemo
 
One of the many consequences of this theorem of Magnus is the following
characterisation of finitely generated free groups, which
Magnus proved a few years later in \cite{13}
\proclaim{Theorem 2}
If the group G can be be generated by q elements, where q is finite,
and if $G/\gamma_n(G) \cong F/\gamma_n(F)$,
where F is a free group of rank q, then $G \cong F$. 
\endproclaim
Theorem 2 will play a role here in due course.
\bigskip
\head 2. One-relator groups
\endhead 
My next encounter with Magnus came some years later. After spending
a postdoctoral year in Manchester, as a 
Special Lecturer, I came to Princeton in 1959 as
an instructor. Some time towards the end of that year, Trueman MacHenry,
a doctoral student of Magnus whom I had met in Manchester, came to visit
me in Princeton. He was accompanied by Bruce Chandler, another of Magnus'
doctoral students. MacHenry brought greetings from Magnus and an implicit
offer of a job at the Courant Institute. I was very happy at the
prospect of working with Magnus, and came up to New York to give a talk
 early in 1960. This was the first time that I had actually
  met Magnus  and I was very impressed with his quickness and
his vast knowledge, not only of mathematics but of almost everything else
as well.
He must have been amused by my talk, which was to an audience
which consisted  mainly of analysts, because I talked about an extremely
esoteric theorem that Norman Blackburn and I \cite {4} had proved about 9 months
earlier, a theorem that only the most special of specialists would have
found interesting. Nonetheless,  the analysts were apparently
convinced by Magnus that
I should be offered a position and the offer was made explicit soon after.
I immediately accepted the position and came to the Courant in the late
summer of 1960. I spent part of that summer in Pasadena, where I met Danny
Gorenstein and Roger Lyndon. Roger Lyndon had been working on 
groups with parametric exponents \cite{9}. There were other things
on his mind as well  and he asked me a number of questions
about one-relator groups while we were in Pasadena. I had not read Magnus's 
papers on one-relator groups, and so I was
forced to spend part
of that summer trying to understand Magnus's work. 
 There are two main theorems that I would like to  describe. The first of these
is Magnus' celebrated "Freiheitssatz"  \cite{10}.
\proclaim{Theorem 3}
Let
$$G = <x_1,\dots,x_q;r=1>$$
be a group defined by a single relator r. If the first and last
letters of r  are not inverses and if $x_1$ appears in r, then
the subgroup of G generated by $x_2,\dots,x_q$ is a free group, freely
generated by $x_2,\dots,x_q$.
\endproclaim
One of the byproducts of the proof of the Freheitssatz was an extraordinary
unravelling of the structure of these groups, which allowed Magnus to
deduce, in due course, that one-relator groups have solvable word problem
\cite{11}:
\proclaim{Theorem 4}
Let G be a group defined by a single relation. Then G has a solvable
word problem.
\endproclaim
\head 3. Surface groups \endhead
As I remarked earlier, I came to the Courant Institute 
in the late summer of 1960.  Magnus was part of the electromagnetic
group of Morris Kline and used to go to a seminar on electromagnetism.
I remember one occasion, after tea, when Morris Kline and Magnus and I were in
the elevator, in the old hat factory which had become
 the Courant Institute. Kline was talking to Magnus and addressed him as
Bill. I saw Magnus shudder at the prospect of being called Bill. Magnus
was much too polite to say anything, but I did not feel at all constrained
to be silent. So I said something like this to Kline: "Morris, Wilhelm
is a Geheimrat and so he simply cannot be called Bill. Wilhelm is more
appropriate". Morris Kline smiled and indicated that this was okay
with him. To which Magnus responded by muttering under his breath:
"Thank God". He later told me that he had hated being addressed as Bill,
and was grateful to me for telling Kline to call him Wilhelm.
Magnus had a large group of Ph.D. students in 1961. They included
Karen Fredericks, Bruce Chandler, Seymour Lipschutz and  Trueman MacHenry,
all of whom worked with Magnus on Combinatorial Group Theory. 
Martin Greendlinger
had already completed his beautiful work on small cancellation theory,
but was still around. In addition, Magnus had some other students who
worked with him on Hill's equation. There was always a line of students
outside his office and I offered to take on some of them to relieve
him of some of
the burden. Magnus ran a seminar on Combinatorial Group Theory at
the time. The seminar was a mixture of pure research and reports by
students  on papers that they were reading. Shortly after I arrived,
the research part of the seminar was broadened and the role of students was
reduced essentially to zero. One of Magnus' habits was to propose a 
number of problems in the seminar. He was always
interested in purely algebraic proofs of theorems that had been 
proved by other means. In particular, he asked late in 1961, whether
there was a direct proof of the residual finiteness of the fundamental
groups of two-dimensional orientable surfaces. Both Karen Frederick 
and I began to work on this problem and we both eventually,
independently,
 came up with
a solution. Her solution \cite{7} was very closely tied to the actual
presentation of these surface groups.  I took a somewhat more general approach
which yielded somewhat more \cite{2}.
The net result was the following
\proclaim{Theorem 5}
Let F be a free group on X, A a free abelian group on Y,
 f an element of F which
is not a proper power in F, and a an element of A which is not a proper
 power in A.
Then the  group
$$G = <X \cup Y; [y,y^\prime]=1 (y,y^\prime \in Y), f=a>$$
is residually free.
\endproclaim
 In particular, it follows from this theorem that surface groups
are residually free and hence residually torsion-free nilpotent (and also
residually finite). 
 The residual nilpotence of one-relator groups
was itself a topic of some interest to Magnus. He later put Bruce Chandler
to work on this topic. Chandler \cite{6} eventually found an alternative
proof of the residual  torsion-freeness of surface groups by making use
of the ring R that I described at the outset. Some years later I found
yet another means of proving the residual nilpotence of surface groups.
I never discussed this method with Magnus, and subsequently forgot about it.
A few months ago, J. Lewin told me that he had also found
a very simple argument
to prove the residual nilpotence of surface groups. This prompted me to
rethink some  of my old ideas and it is to these thoughts that I want
to turn  next.
\head 4. Residually nilpotent groups
\endhead
The fundamental groups of two dimensional orientable surfaces all contain
a free subgroup with infinite cyclic factor group. Thus they have the
same form as the groups covered by the following theorem.
\proclaim{Theorem 6} Let G be a finitely generated group. Suppose that
G contains a free, normal subgroup N such that G/N is infinite cyclic.
If $G/\gamma_2(N)$ is residually torsion-free nilpotent, then so is G.
\endproclaim
This theorem is similar in spirit to a theorem of P. Hall \cite{8}, who 
proved that if $G$ is a group with a normal, nilpotent subgroup $N$,
then $G$ is nilpotent if $G/\gamma_2(N)$ is nilpotent.
 Theorem 6 leads to a host of new examples of residually nilpotent
groups. In particular, it can  be used to give yet another proof
of the residual torsion-free nilpotence of surface groups.
 There is a further use of Theorem 6 that I want to describe here.
To this end, let me recall the following definition.
\definition{Definition} A group G is termed parafree if G is residually
nilpotent and there exists a free group F such that
$$G/\gamma_n(G) \cong F/\gamma_n(F) \ for \ all\ n.$$
\enddefinition
 
Parafree groups, which can be likened to free groups, exist in 
profusion, see, e.g., \cite{3}. It should be pointed out that it
 follows from  Magnus' Theorem 2 that a non-free, finitely generated
parafree group $G$ with the same nilpotent factor groups $G/\gamma_n(G)$
as a free group of rank $q$ cannot be generated by $q$ elements.
It is not known how closely a parafree group can resemble a free group.
I want to describe next some new, non-free parafree groups, which
very closely resemble free groups. The proof that these groups
are parafree is an easy application of Theorem 5 (see {\bf 6}).
\proclaim{Theorem 7} Let F be the free group on 
$s,t,a_1,\dots,a_q$ and let $w$ be an element of F which involves
$a_1$ and does not
involve $s$. In addition suppose that $w$ lies in the $k-th$ term
$F^{(k)}$ of 
the derived series of $F$. Then the one-relator group
$$G_w = < s,t,a_1,\dots,a_q;a_1=ws^{-1}t^{-1}st>$$
is parafree and not free. Moreover
$$G/G_w^{(k)} \cong H/H^{(k)},$$
where H is a free group of rank $q+1$.
\endproclaim
I will sketch the proofs of Theorems 6 and 7 in section {\bf 6}.
\bigskip
\head 5. Some problems on $\Cal D$-groups 
\endhead
 Let $p_1=2,\ p_2=3,\dots $ be the set of all primes in ascending order
of magnitude.
\definition{Definition} A group G is called a $\Cal D$-group if it admits a set of unary operators
$$\pi_1,\pi_2,\dots$$
such that for all  $g \in G$
$$g^{p_i}\pi_i=g=(g\pi_i)^{\pi_i}.$$
\enddefinition
It is not hard to verify that these $\Cal D$-groups consist precisely of
those groups $G$ in which extraction of $n$-th roots is uniquely
possible, for every positive integer $n$.
The class of all $\Cal D$-groups form a {\it variety of (universal) algebras}. 
The precise technical description of these terms does not matter. It
suffices only to say that in such a variety one has the notion of a free $\Cal D$-group,
as well as all of  the other notions that one makes use of in group theory.
 There are a number of properties of such free $\Cal D$-groups that
are similar to properties of   free groups. For example one has the following
\proclaim{Theorem 8}
$\Cal D$-subgroups of free $\Cal D$-groups are free.
\endproclaim
I was told about this theorem by Tekla Taylor-Lewin, but I have not
been able to locate a reference.
Magnus was fond of these $\Cal D$-groups, and discussed them in his
book with Karrass and Solitar \cite{14}.
It seems appropriate, therefore, to raise here  some new problems
 about free $\Cal D$-groups,
in his  memory.
  To this end, let $F$ be the free $\Cal D$-group on
$$x_1,\dots,x_q.$$
It is not hard to see that the mapping
$$x_i \mapsto 1+\xi_i \ (i=1,\dots,q)$$
defines a homomorphism $\phi$ of $F$ into the group of units of $R$.
\proclaim{Problem 1} Is $\phi$ a monomorphism?
\endproclaim
\proclaim{Problem 2} Let 
$$G=<x_1,\dots,x_q;r=1>.$$
If extraction of n-th roots in G is unique, whenever such roots exist,
can G be embedded in a $\Cal D$-group?
\endproclaim
\proclaim{Problem 3} Suppose that the one-relator group 
$$G=<x_1,\dots,x_q;r=1>$$
can be embedded in a $\Cal D$-group. If H is the one-relator
$\Cal D$-group generated, as a $\Cal D$-group, by $x_1,\dots,x_q$
and defined, as a $\Cal D$-group by the single relation $r=1$,
is the word problem solvable for H? In general, is the word 
problem solvable for one-relator $\Cal D$-groups?
\endproclaim
\proclaim{Problem 4} Is there a freiheitssatz for one-relator 
$\Cal D$-groups?
\endproclaim
\proclaim{Problem 5} Can free groups be characterised by a length
function?
\endproclaim
\proclaim{Problem 6} Does a free $\Cal D$-group act freely on a
$\Lambda$-tree, for a suitable choice of $\Lambda$?
\endproclaim                                         
\head{6. Proofs}
\endhead
I want to sketch here the proofs of Theorem 6 and Theorem 7.
I would like to begin with the proof of Theorem 6. I will
adopt the notation used in the statement of the theorem. Since
$G/N$ is infinite cyclic, we can choose an element $t \in G$, such
that $G=gp(N,t)$. So $t$ is of infinite order modulo $N$. Suppose
that $g \in G, g\neq 1$. We want to find a normal subgroup
$K$ of $G$  such that $G/K$ is torsion-free nilpotent and $g \notin
K$. It is clear that it suffices to consider the case where
$g \in N$. Since $N$ is free there exists an integer $k$ such that
$g \notin \gamma_k(G)$. Notice that $H=N/\gamma_k(N)$ is a torsion-free
nilpotent group and so we can form the Malcev completion
$\overline{H}$ of $H$ (see \cite{15}). This group $\overline{H}$ is
a minimal $\Cal D$-group containing $H$. It is again nilpotent of
class $k-1$
and torsion-free. Moreover if we denote by $\tau$
 the automorphism that $t$ induces on $\overline{H}$, then
$\tau$ extends uniquely to an automorphism $\overline{\tau}$
of $\overline{H}$. Let $\overline{G}$ be the semidirect product
of $\overline{H}$ with the infinite cyclic group generated by
an element $\overline{t}$, where $\overline t$ acts on $\overline H$
as $\overline{\tau}$. Then $M=\overline H /\gamma_2(\overline H)$ is
a direct sum of copies of the additive group of $\Bbb Q$ (see
\cite{1}). Morevover, it is not hard to deduce from the fact that
$G/\gamma_2(N)$
is residually torsion-free nilpotent, that
 $\overline G/\gamma_2(\overline H)$ is also residually torsion-free
nilpotent. We now view $M$  as
a module over the rational group algebra $\Gamma$ of the infinite
cyclic group on $\overline t$. Since $\Gamma$ is a principal ideal
domain and the group $G$ is finitely generated, the $\Gamma$-module
$M$ is finitely generated and hence a direct sum of cyclic modules.
This decomposition of $M$ makes it possible to understand the
action of $\overline t$ on $\overline H$ and, using the fact that
$H=F/\gamma_k(F)$, to deduce that $\overline H$ is residually
torsion-free nilpotent. It is easy to deduce that $G$ is itself
residually torsion-free nilpotent.
There are three steps in the proof of Theorem 7. 
The first makes use of Magnus'
method of unravelling the structure of one-relator groups. 
This method allows to prove, easily, that 
the normal closure $N$ of $s,a_1,\dots,a_q$ in $G_w$ is free.
The second step involves the verification that  $G/\gamma_2(N)$ is residually
 torsion-free nilpotent. So, by Theorem 6, $G_w$ is residually
torsion-free nilpotent. The other properties of $G_w$ follow directly
from the form of the defining relation of $G_w$. The final step in the
proof of Theorem 7 is the verification that $G_w$ is not free. This is
accomplished by invoking an algorithm of J.H.C. Whitehead \cite{15}.
I have only talked here about a few of Magnus' theorems. All of
his work is
filled with beautiful, new ideas, giving joy to all of us. 
Wilhelm Magnus is sorely missed, but his work will be with us always.

\enddocument
\end